\documentclass[english]{amsart}

\usepackage{babel}
\usepackage{amstext}
\usepackage{amsmath}
\usepackage{amsfonts}
\usepackage{latexsym}
\usepackage{ifthen}
\usepackage{xypic}
\xyoption{all}
\newcommand{\E}{{\mathcal E}}

\newcommand{\sJ}{{\mathcal J}}
\newcommand\sE{{\mathcal E}}
\newcommand\sME{{{\overline {ME}}}}
\newcommand\sF{{\mathcal F}}

\newcommand\sI{{\mathcal I}}

\newcommand\sL{{\mathcal L}}
\newcommand\sO{{\mathcal O}}
\newcommand\sS{{\mathcal S}}

\newcommand\sC{{\mathcal C}}
\newcommand\sT{{\mathcal T}}
\newcommand\bR{{\mathbb R}}

\newcommand\bQ{{\mathbb Q}}

\newcommand\bP{{\mathbb P}}

\newcounter{lemma}

\newtheorem{lemma1}[lemma]{\setcounter{equation}{0}}

\newenvironment{lemma}{\begin{lemma1}{\bf Lemma.}}{\end{lemma1}}

\newenvironment{example}{\begin{lemma1}{\bf Example.}\rm}{\end{lemma1}}

\newenvironment{theorem}{\begin{lemma1}{\bf Theorem.}}{\end{lemma1}}
\newenvironment{question}{\begin{lemma1}{\bf Question.}}{\end{lemma1}}

\newenvironment{proposition}{\begin{lemma1}{\bf Proposition.}}{\end{lemma1}}

\newenvironment{corollary}{\begin{lemma1}{\bf Corollary.}}{\end{lemma1}}

\newenvironment{remark}{\begin{lemma1}{\bf Remark.}\rm}{\end{lemma1}}
\newenvironment{definition}{\begin{lemma1}{\bf Definition.}}{\end{lemma1}}

\newenvironment{notation}{\begin{lemma1}{\bf Notation.}}{\end{lemma1}}

\newenvironment{conjecture}{\begin{lemma1}{\bf Conjecture.}}{\end{lemma1}}

\newenvironment{Induction Step}{\begin{lemma1}{\bf Induction Step.}}{\end{lemma1}}
\newenvironment{Proof of Theorem 1.2}{\begin{lemma1}{\bf Proof of Theorem 1.2.}}{\end{lemma1}}

\begin{document}

\title {Varieties with generically nef tangent bundles}

\author[Th. Peternell]{Thomas Peternell}
\address{Th. Peternell - 
Mathematisches Institut - Universit\"at Bayreuth - D-95440 Bayreuth, Germany} 
\email{thomas.peternell@uni-bayreuth.de}
\date{\today}
\maketitle
\tableofcontents

\section{Introduction}
Generic nefness is a concept which appeared first in Miyaoka's important paper [Mi87], in which he studied the cotangent bundles of non-uniruled manifolds: 
\vskip .2cm \noindent 
{\it the cotangent bundle of a non-uniruled manifold is generically nef}. 
\vskip .2cm \noindent 
A vector bundle $E$ on a projective $n-$dimensional manifold $X$ is {\it generically nef}, if the following holds. Given any ample line bundles $H_j, 1 \leq j \leq n-1,$
let $C$ be a curve cut out by general elements in $\vert m_j H_j \vert$ for $m_j \gg 0,$ then the restriction $E \vert C $ is nef. \\
In the same way {\it generic ampleness} is defined. If we fix the $H_j$, we speak of generic nefness with respect to $(H_1, \ldots, H_{n-1}).$ 
\vskip .2cm \noindent
In this paper we are mainly interested in manifolds with generically nef {\it tangent} bundles, and the aim is to better understand uniruled varieties, in particular
those with nef anti-canonical bundles. \\
A projective manifold $X$ whose tangent bundle $T_X$ generically nef w.r.t some polarization is uniruled unless $K_X \equiv 0$; this is a special
case of Miyaoka's theorem.
Actually much more holds: 
if $C \subset X$ is an irreducible curve with
$K_X \cdot C < 0$ and $T_X \vert C$ is nef, then $X$ is uniruled. If $T_X \vert C$ is even ample, then $X$ is rationally connected. Uniruledness however is by far not
enough to have a generically nef tangent bundle. This is demonstrated by the following mapping property.
\vskip .2cm \noindent 
{\it If $T_X$ is generically nef with respect to some polarisation and if $f: X \to Y$ is a surjective holomorphic map to a normal projective variety $Y,$ then $Y$ is uniruled
or a smooth model of $Y$ has Kodaira dimension $0.$ }  
\vskip .2cm \noindent
\vskip .2cm \noindent
On the other hand there is of course no restriction for images of uniruled varieties.
We are then concerned with 
the question to which extent the converse of the mapping property holds. Before we address that problem, we discuss another nefness concept for vector bundles.
The reason for doing so is that especially for manifolds with nef anti-canonical bundles the notion of generic nefness is not good enough, 
because it
gives in some sense only birational information. Therefore we introduce the notion of a sufficiently nef vector bundle $E$: 
\vskip .2cm \noindent a vector bundle $E$ is {\it sufficiently nef} if through any point of $X$ there is
a family of (generically irreducible) curves covering the whole manifold $X$ such that the bundle $E$ is nef on the general member of this family. \\
In the same way
we define {\it sufficiently ample} bundles. 
\vskip .2cm \noindent
There is a relation between both ampleness notions:
\begin{theorem} If $E$ is generically ample w.r.t some polarization, then $E$ is sufficiently ample. 
\end{theorem} 
This is proved in section 7. Unfortunately the analogous statement in the nef case is false for general bundles. \\
The importance of the notion of sufficient nefness lies in particular in 
\begin{theorem} If $T_X$ is sufficiently nef, then the Albanese map is a (surjective) submersion. 
\end{theorem} 
It is conjectured that the Albanese map of a projective manifold $X$ with $-K_X$ nef is a (surjective) submersion, the surjectivity being known by
Qi Zhang [Zh05]. Moreover manifolds with nef anti-canonical bundles have the above stated mapping property. Hence we are lead to the
following
\begin{conjecture} 
Let $X$ be a projective manifold. 
\begin{enumerate} 
\item If $-K_X$ is nef, then $T_X$ is generically nef.
\item If $-K_X$ is nef, then $T_X$ is sufficiently nef.
\item If $X$ is rationally connected, then $T_X$ is generically ample (hence also sufficiently ample).
\end{enumerate} 
\end{conjecture}

Theorem 1.2 says that (1.3.2) implies the smoothness of the Albanese map of a manifold with nef anti-canonical class. 
In this paper we prove Conjecture 1.3 in the following cases.
\begin{theorem}
\begin{itemize}
\item (1.3.1) and (1.3.2) hold if $-K_X$ is semi-ample, i.e., some multiple $-mK_X$ is generated by global sections. 
\item (1.3.3) holds if $-K_X$ is big and nef, in particular for Fano manifolds $X.$
\item (1.3.3) holds for rational surfaces and suitable ample $H.$ 
\end{itemize} 
\end{theorem} 

If $X$ is rationally connected, we prove in all dimension that there is an irreducible $C$ whose numerical equivalence class is in the interior 
of the movable cone, such that $T_X \vert C$ is ample. This can be seem as a weak substitute for (1.3.2). 

\vskip .2cm \noindent The verification of the conjecture for
Fano manifolds with $b_2 > 1$ uses a cone theorem for the movable cone proved by Birkar-Cascini-Hacon-McKernan [BCHM06] and further results from that
paper.

\vskip .2cm \noindent An important issue is the relation to stability. If a vector bundle $E$ is semi-stable w.r.t. $h = (H_1, \ldots, H_{n-1})$
and if $\det E \cdot h \geq 0,$ 
then it is easy to see that $E$ is generically nef w.r.t. $h.$ Thus the generic ampleness of the tangent bundle of a Fano manifold $X$ with $b_2 = 1$
can be seen as a weak substitute for the potential stability of $T_X$. Moreover, we use methods from the theory of semi-stable vector
bundles, namely the Mehta-Ramanathan theorem, to prove generic ampleness of the tangent bundle of a Fano manaifold.

\vskip .2cm \noindent As a consequence of the above results we conclude that the Albanese map for a manifold with semi-ample anti-canonical bundle
is a surjective submersion, that all section in a tensor of the cotangent bundle do not allow zeroes and that the rational quotient is a submersion
at all points where it is defined. These properties also follow from [DPS96], since the anti-canonical bundle is hermitian semi-positive 
(using K\"ahler geometry), but the method presented here has the advantage to immediately generalize also to singular situations. For simplicity we
will however formulate all results only in the smooth case.   \\
For important discussions concerning parts of section 6 I would like to thank James McKernan. 

\section{Basic definitions}
\setcounter{lemma}{0}

Unless otherwise stated, $X$ always denotes a projective manifold of dimension $n.$ 

\begin{definition} {\rm 
\begin{enumerate} 
\item Let $H_1, \ldots, H_{n-1} $ be ample divisors. 
A vector bundle $E$ is said to be $(H_1, \ldots, H_{n-1}) -$ {\it generically nef (ample)}, if $E \vert C$ is
nef (ample) for a general curve $C = D_1 \cap \ldots \cap D_{n-1}$ for general $D_i \in \vert m_i H_i \vert $ and $m_i \gg 0.$ 
Such a curve is called {\it MR-general}, which is to say ``general in the sense of Mehta-Ramanathan''. 
\item The vector bundle $E$ is called {\it generically nef (ample)}, if $E$ is  $(H_1, \ldots, H_{n-1}) -$ generically nef (ample) for all $H_i.$
\item  $E$ is {\it almost nef [DPS00]}, if there is a countable union $S$ of algebraic subvarieties such $E \vert C$ is nef for all curves
$C \not \subset S.$ 
\end{enumerate} }
\end{definition}

If $Y$ is a normal projective variety and $L$ a Weil divisor on $Y$, then $L$ is locally free on the regular part of $Y$. Let $C$ be MR-general w.r.t 
$(H_1, \ldots, H_{n-1}).$ 
Then $C$ does not meet the singular locus of $Y$ and thus $L \cdot C $ makes perfect sense. {\it Hence the notion of generic nefness (ampleness) is 
well-defined
also for Weil divisors}.  

\begin{notation} {\rm Let $X$ be a normal variety with smooth locus $X_0$ and inclusion $i: X_0 \to X.$ Then we define the tangent sheaf
$$ \sT = i_*(\sT_{X_0}) $$
and the sheaf of $p-$forms 
$$ \tilde \Omega^p_X = i_*(\Omega^p_{X_0}). $$ } 
\end{notation} 

\begin{notation} {\rm Fix ample line bundles $H_i$ on a projective variety $X$ and let $E$ be a vector bundle on $X$. Then we define the slope
$$ \mu_{(H_1, \ldots, H_{n-1})}(E) = \det E \cdot H_1 \ldots \cdot H_{n-1} $$
and obtain the usual notion of (semi-)stability w.r.t. $(H_1, \ldots, H_{n-1}).$ }
\end{notation} 

The importance of the notion of MR-generality comes from Mehta-Ranamathan's theorem [MR82] (see Proposition 5.1): 

\vskip .2cm \noindent {\it If $E$ is semi-stable w.r.t. $(H_1, \ldots, H_{n-1})$ and if 
$c_1(E) \cdot H_1, \ldots, H_{n-1}) \geq 0,$ then $E$ is generically nef w.r.t. $(H_1, \ldots, H_{n-1})$. }
\vskip .2cm \noindent


\begin{definition} A vector bundle $E$ is numerically flat [DPS94] iff $E$ is nef with $c_1(E) = 0.$ 
\end{definition} 

Numerically flat bundles are filtered by hermitian flat bundles (which by definition are given by unitary representations of the fundamental group):

\begin{theorem} A vector bundle $E$ is numerically flat iff there is a filtration
$$ 0 = E_0 \subset E_1 \subset \ldots \subset E_p = E $$
of subbundles such that the graded pieces $E_{j+1}/E_{j} $ are hermitian flat.
In particular $E$ is semi-stable w.r.t any polarization. 
\end{theorem}

For the proof we refer to [DPS94]. 

\begin{notation} {\rm As usual, $\kappa (L)$ is the Kodaira dimension of the line bundle $L$, and $\nu(L)$ denotes the numerical dimension
of a nef line bundle $L,$ i.e., the largest number $m$ so that $L \not \equiv 0.$ 
The Kodaira dimension of a normal projective variety will be the Kodaira dimension of a desingularisation. \\
Given a projective manifold $X$, the movable cone is denoted by $\sME(X),$ see [BDPP04]. Also the subcone 
$ \sME(X)^{-}$ consisting of those $\gamma \in \sME(X)$ for which $K_X \cdot \gamma \leq 0$ will be of importance. \
The closed cone generated by the MR-general curves will be denoted by $\mathcal{CI}(X)$ (``complete intersection''). 
}
\end{notation}

\section{Generically nef divisors}
\setcounter{lemma}{0}

In this section we study varieties with generically nef anti-canonical divisors. 

\begin{theorem} Let $Y$ be a normal projective variety. Assume that $-K_Y$ is generically nef w.r.t. $(H_1, \ldots, H_{n-1}).$   
Then either $Y$ is uniruled or $\kappa (Y) = 0.$ 
\end{theorem} 

\begin{proof} Let $\pi: \hat Y \to Y$ be a desingularisation. Let $C \subset Y$ be MR-general  w.r.t. $(H_1, \ldots, H_{n-1})$   
and let $\hat C$ denote its strict transform in $\hat Y.$
Since $C$ does not meet the center of $\pi,$ we have
$$ K_{\hat Y} \cdot \hat C = K_Y \cdot C \leq 0.$$
If $K_{\hat Y} \cdot \hat C < 0$, then $\hat Y$ (hence $Y$) is uniruled by [MM86]. Therefore
we may assume that
$$ K_{\hat Y} \cdot \pi^*(H_1) \cdot \ldots \cdot \pi^*(H_{n-1}) = 0. \eqno (*) $$
We assume that $\hat Y$ is not uniruled and must show $\kappa (\hat Y) = 0.$ So by [BDPP04] $K_{\hat Y}$ is pseudo-effective. 
Hence $c_1(K_{\hat Y})$ is represented by a positive closed current $T$, e.g. via a singular metric on $K_{\hat Y}.$ 
From (*) we deduce
$$ \pi_*(T) = 0.$$ 
Therefore the support of $T$ is contained in the exceptional locus of $\pi$, hence by Siu's theorem [Si74], we obtain a
decomposition
$$ T = \sum a_i T_{E_i}, $$
where $a_i > 0,$ the $E_i$ are $\pi-$exceptional and $T_{E_i}$ denotes the current ``integration over $E_i$''.
By [BDPP04,3.7,3.10], the $a_i$ are rational, so that some multiple $mK_{\hat Y} $ is numerically equivalent to an effective divisor. 
From [CP07] finally, we obtain $\kappa (\hat Y) \geq 0$, hence $\kappa (\hat Y) = 0.$ 
\end{proof}

An application of Theorem 3.1 is 

\begin{corollary} Let $X$ be a projective manifold, $f: X \to Y$ a surjective holomorphic map to a normal projective variety. Assume one of the
following.
\begin{enumerate}
\item $-K_X$ is nef;
\item $-mK_X = \sO_X(D)$ for some positive $m$ and $(X,D)$ is log canonical;
\item there exists a positive number $m$ such that $h^0(-mK_X) > 0$ and the base locus of $\vert -mK_X \vert $ does not project onto $Y$.
\end{enumerate}
Then either $Y$ is uniruled or $\kappa (Y) = 0.$ In particular the Albanese map of $X$ is surjective. 
\end{corollary} 

\begin{proof} By [Zh05], $-K_Y$ is generically nef. This is not explicitly stated in [Zh05], but this is what the proofs of Theorems 2,3 and 4 in [Zh05]
actually give. Now apply 3.1.
\end{proof} 

In [Zh05], even more is shown in case $-K_X$ nef: if $f: X \dasharrow  Y$ is a dominant rational map to the projective manifold $Y$, then $Y$ is uniruled
or $\kappa (Y) = 0.$

It is natural to ask for a common generalization for the second and third assumptions in Theorem 3.2; if $f$ is a submersion, this is contained in 
[DPS01,theorem 2.11]; if $Y$ is smooth and $f$ is flat, it follows from H\"oring [Hoe06]. We generalize this as follows.

\begin{theorem} Let $f: X \to Y$ be a surjective holomorphic map from the projective manifold to the normal projective variety $Y.$
Assume that $-K_X$ is pseudo-effective and let $h$ be a singular metric on $-K_X$ with positive curvature current. Assume that the support
of the multiplier ideal $\sI(h)$ does not project on $Y$. Then $-K_Y$ is generically nef. 
\end{theorem} 

\begin{proof} Let $\pi: \hat Y \to Y$ be a desingularisation; choose a birational map $\pi: \hat X \to X$ with $\hat X$ such that the induced
map $\hat f: \hat X \dasharrow \hat Y$ is holomorphic. Then the main result of Berndtsson-Paun [BP07] implies that 
$$ K_{\hat X/\hat Y} + \pi^*(-K_X) $$
is pseudo-effective. Write $K_{\hat X} = \pi^*(K_X) + \sum a_i E_i$ with $a_i > 0$ and $E_i$ being $\pi-$exceptional. Then we conclude that
$$ \hat f^*(-K_{\hat Y}) + \sum a_i E_i$$
is pseudo-effective. This clearly gives our claim by (3.1).  
\end{proof} 

Notice that the support condition is really necessary; in fact, Zhang observed that there are many ruled surfaces over curves of genus
$\geq 2$ with big (but not nef) anti-canonical bundle. 

\begin{question} Let $f: X \to Y$ be a surjective map of projective manifolds. Assume $-K_X$ nef. Is $-K_Y$ pseudo-effective or even nef?
\end{question} 

In general, a generically nef line bundle is not pseudo-effective. The reason is simply that the movable cone $\sME(X)$ in general is larger than the
closed cone $\sC \sI(X)$ generated by MR-general curves. See [BDPP04] for an example. However we might ask:

\begin{question} Suppose $X$ is a projective manifold with $K_X$ generically nef. Is $K_X$ pseudo-effective? In other words, if $K_X$ is not pseudo-effective,
are there ample divisors $H_i,$ such that
$$ K_X \cdot H_1 \cdot \ldots \cdot H_{n-1} < 0? $$
\end{question}  

See [CP98] for some positive results. 

\section{Nefness properties of the tangent bundle}
\setcounter{lemma}{0}

Miyaoka-Mori [MM86] proved that if $K_X \cdot C < 0 $ for some MR-general curve (actually a curve moving in a covering family suffices), then $X$ is uniruled. This can be seen as a special case
of Miyaoka's theorem [Mi87] stating that if $T^*_X \vert C$ is not nef on the MR-general curve $C = m_1H_1 \cap \ldots \cap m_{n-1}H_{n-1},$ then 
$X$ is uniruled. It is natural to ask whether one really needs to consider MR-general curves:

\begin{question} Let $T \subset {\rm Chow}(X)$ be an irreducible component and assume that $T^*_X \vert C_t$ is not nef for general $t \in T.$
Is $X$ uniruled?
\end{question} 

\begin{remark} {\rm It is not sufficient to require $T^*_X \vert C_t$ to be nef for the general member $C_t$ of some covering 
family of curves. In fact, in
[BDPP04] it is shown that any K3 surface and any Calabi-Yau threefold admit a covering family $(C_t)$ of curves such that $T^*_X \vert C_t$ is not
nef for general $t.$ In other words $T^*_X$ (as well as $T_X$) is not almost nef. Consequently we also find projective $n-$folds $X$ with 
$\kappa (X) = n-2$ such that $T^*_X$ is not almost nef. }
\end{remark} 

\begin{question} Let $X$ be a projective manifold with $\kappa (X) \geq \dim X - 1$. Is $T^*_X$ almost nef? 
\end{question}

Since Question 4.1 seems to be too hard at the moment, we consider the stronger condition $T_X \vert C$ to be nef. In [Pe06] the following is proved, for the ampleness
statement see also [KST07],[BM01]. 

\begin{theorem} Let $X$ be a projective manifold, $C \subset X$ an irreducible curve. If $T_X \vert C$ is nef, then $\kappa (X) < dim X.$ 
If additionally $K_X \cdot C < 0,$ then $X$ is uniruled. If $T_X \vert C$ is ample, then $X$ is rationally connected.
\end{theorem} 

\begin{corollary} Suppose $X$ uniruled with rational quotient $f: X \dasharrow W,$ where $W$ is smooth. Let $C \subset X$ be an irreducible curve such that
$f$ is holomorphic near $C$ and assume $\dim f(C) > 0.$ Suppose that $T_X \vert C$ is nef and let $\sS \subset T_X \vert C$ be the maximal ample subsheaf [PS00]. 
If ${\rm rk}\sS = m > 0,$ then $\dim W \leq n-m.$
\end{corollary}

\begin{proof} We simply observe that the natural morphism
$$ \sS \to f^*T_W \vert C $$
must vanish, because otherwise $T_W \vert f(C)$ would contain an ample subsheaf. Since $T_W \vert f(C)$ is nef (the map $T_X \vert C \to f^*T_W \vert C$
being generically surjective), we conclude by 4.4 that $W$ is uniruled, a contradiction. 
\end{proof} 

The next theorem shows that manifolds with generically nef tangent bundles have the same mapping properties as manifolds with nef 
anti-canonical bundles. 
\begin{theorem} Assume that $T_X$ is generically nef w.r.t $(H_1, \ldots, H_{n-1})$. Let $f: X \to Y$ be a surjective holomorphic map to the normal
projective variety $Y$. 
Then $Y$ is uniruled or $\kappa (Y) = 0.$
\end{theorem} 

\begin{proof} Let $C \subset X$ be MR-general. Fix a desingularisation $\tau: \hat Y \to Y$ and choose a smooth birational model $\pi: \hat X \to X$ 
such that the induced map $\hat f: \hat X \to \hat Y$ is holomorphic. Then $C$ can be considered as curve on $\hat X$ (but it is not a MR-curve on $\hat X$), and
$T_{\hat X} \vert C$ is nef.
We claim that
$$ K_{\hat Y} \cdot \hat f(C) \leq 0. \eqno (*)$$
In fact, otherwise we find a positive integer $N$ and a non-zero section
$$ \omega \in H^0(NK_{\hat Y} \vert \hat f(C)) = H^0((\Omega^m_{\hat Y})^{\otimes N} \vert \hat f(C))$$ 
having a zero ($ m = \dim \hat Y)$. Hence we obtain a section
$$ \hat f^*(\omega) \in H^0((\Omega^m_{\hat X})^{\otimes N} \vert C) $$
with a zero. This contradicts the nefness of $(\bigwedge^mT_{\hat X})^{\otimes N} \vert C.$ \\
Now we proceed similarly as in (3.1). If we have strict inequality in (*), then $Y$ is uniruled. So suppose
$$ K_{\hat Y} \cdot \hat f(C) = 0.$$ 
Assuming $Y$ not uniruled, we must show that $\kappa (\hat Y) = 0.$ By [BDPP04], $K_{\hat Y}$ is pseudo-effective. 
Hence $L = \hat f^*(K_{\hat Y})$ is pseudo-effective with $L \cdot C = 0.$ Then we just follow the arguments of (3.1) to conclude
that $\kappa (\hat Y) = 0.$

\end{proof}

In particular the Albanese map of a manifold $X$ whose tangent bundle $T_X$ is generically $(H_1, \ldots, H_{n-1})$-nef is surjective. Now we ask for conditions
under which the Albanese map is a submersion. Here we naturally need conditions concerning curves through every point. 

\begin{definition} We say that $T_X$ is sufficiently nef (ample), if the following holds. Given any point $x \in X,$ there exists a covering family $(C_t)$ of (generically irreducible) curves
through $x$ such that $T_X \vert C_t$ is nef (ample) for general $t.$ 
\end{definition} 

In this notation we have 

\begin{proposition} If $T_X$ is sufficiently nef, then the Albanese map of $X$ is a surjective submersion. 
\end{proposition} 

\begin{proof} If the Albanese map is not surjective or not a submersion, then we find a nonzero holomorphic 1-form $\omega$ on $X$ admitting a zero $x.$ 
Now consider a curve $C$ through $x$ such that $T_X \vert C$ is nef and such that $\omega \vert C \in H^0(\Omega^1_X \vert C)$ does not vanish identically. Such a curve
exists by our assumption. Since $\omega(x) = 0,$ we obtain a contradiction.
\end{proof} 

\begin{corollary} If $T_X$ is almost nef, then the Albanese is a surjective submersion.
\end{corollary} 

\begin{remark} {\rm  The arguments of (4.8) show more generally the following. If 
$$0 \ne \omega \in H^0(X,(\Omega^1_X)^{\otimes N}),$$
then $\omega $ has no zeroes. As an application, if $f: X \to Y$ is a surjective holomorphic map to a projective manifold $Y$ with $\kappa (Y) \geq  0$
(where $T_X$ is sufficiently nef)
then $mK_Y = \sO_Y$ for some $m > 0$ (and
any non-zero $\eta \in H^0(Y,(\Omega^1_Y)^{\otimes N})$ is without zeroes). }
\end{remark} 
Coming back to manifolds with nef anti-canonical bundles, we recall the following

\begin{conjecture} Let $X$ be a projective manifold with $-K_X$ nef. Then the Albanese map is a (surjective) submersion. 
\end{conjecture}

This is known in dimension 3 [PS98] and in all dimensions if $-K_X$ is hermitian semi-positive [DPS96] and, more generally, if $-K_X$ has a singular
metric $h$ with trivial multiplier ideal $\sI(h) = \sO_X$ [DPS01] - even in the K\"ahler case. As already mentioned, the surjectivity is due to Qi Zhang, see [Zh05]. 
However in the K\"ahler case, surjectivity is still unknown. 

Conjecture 4.11 is potentially a consequence of Proposition 4.8 via the first part of the following 

\begin{conjecture} Let $X$ be a projective manifold. 
\begin{enumerate} 
\item If $-K_X$ nef, then $T_X$ is sufficiently nef and also generically nef. \\
\item If $X$ is rationally connected, e.g. Fano, then $T_X$ is generically ample.
\end{enumerate} 
\end{conjecture} 

\begin{example} Let $X$ be a projective surface with $-K_X$ nef. Then $T_X$ is sufficiently nef, as shown by the following
considerations, using classification.  \\
If $q(X) \geq 2,$ then $X$ is a torus, and there is nothing to prove. If $q(X) = 1,$ then either $X$ is hyperelliptic and we are 
done, or $X = \bP(E)$ with a rank 2-vector bundle $E$
over an elliptic curve $A.$ Then automatically $T_X$ is nef [CP91].\\
If $q(X) = 0$ and $X$ is not uniruled, then $X$ is K3 (or Enriques), this being settled by (7.2) below. 
So it remains to treat the rational case. Here we may assume that $X$ is the plane
blown up in at most 9 points $p_1, \ldots, p_r.$ Fix $x \in X.$ If $x$ is disjoint from the exceptional locus of $\pi: X \to \bP_2,$ then the claim
is clear (take lines). Otherwise we have say $\pi(x) = p_1.$ Since $-K_X$ is nef, at most three points can be infinitely near and it follows easily 
that there is a 1-dimensional family of strict transforms of conics through $p_1$ which all contain $x.$  

\end{example}

\section{Semi-stability}
\setcounter{lemma}{0}

In this section we shall bring methods from the theory of stable vector bundles into the game. We start with the following observation. 

\begin{proposition} Let $X$ be a projective manifold of dimension $n$ and $E$ a rank $r-$vector bundle on $X$. 
Suppose that $E$ is $(H_1, \ldots, H_{n-1})-$semi-stable and that
$$\det E \cdot H_1 \cdot \ldots H_{n-1} \geq 0.$$ Then $E$ is  $(H_1, \ldots, H_{n-1})$-generically nef. \\
$E$ is even generically ample w.r.t
$(H_1, \ldots, H_{n-1})$ unless $\det E \cdot H_1 \cdot \ldots \cdot H_{n-1} = 0.$ In that case $E \vert C$ is flat for $C$ MR-general. 

\end{proposition} 

\begin{proof} Let $C$ be MR-general w.r.t. $(H_1, \ldots, H_{n-1}).$ Then $E \vert C$ is semi-stable. By Miyaoka [Mi87] this is equivalent to
saying that the $\bQ-$bundle
$$ (E \otimes {{\det E^*} \over {r}}) \vert C$$
is nef (i.e. flat since the first Chern class of this $\bQ-$bundle vanishes). Since $\det E \cdot C \geq  0,$ we conclude that $E \vert C$ is nef.  
\end{proof} 

\begin{remark} {\rm  Suppose that
$\det E \cdot H_1 \cdot \ldots H_{n-1} = 0.$ Then semi-stability and generic nefness w.r.t $(H_1, \ldots, H_{n-1})$ are the same. 
In fact, if $C$ is MR-general  w.r.t $(H_1, \ldots, H_{n-1})$, so that $E \vert C$ is nef, hence numerically flat ((2.4), possibly normalize) 
and consequently  $E \vert C$ is semi-stable (2.5).
Thus  by [MR82], $E$ is semi-stable w.r.t $(H_1, \ldots, H_{n-1})$. The remaining claim is clear.}
\end{remark} 

\begin{corollary}\begin{enumerate} 
\item  Assume $T_X $ semi-stable w.r.t $(H_1, \ldots, H_{n-1})$. If $-K_X \cdot H_1 \cdot \ldots \cdot H_{n-1} > 0 $, then $X$ is rationally
connected.
\item Assume $T_X$ stable w.r.t $(H_1, \ldots, H_{n-1})$. If $-K_X \cdot H_1 \cdot \ldots \cdot H_{n-1} \geq 0$ and if
$K_X \not \equiv 0 $, then $X$ is rationally connected. 
\end{enumerate}
\end{corollary} 

\begin{proof} (1) If  $-K_X \cdot H_1 \cdot \ldots \cdot H_{n-1} > 0 $, then $T_X \vert C$ is ample for $C$ MR-general, hence we may apply (4.4). \\
(2) By (1) we may assume that $-K_X \cdot H_1 \cdot \ldots \cdot H_{n-1} = 0.$ Since $K_X \not \equiv 0$, we find $H'_1, \ldots, H'_{n-1}$ ample, 
such that for all positive $\epsilon:$
$$ K_X \cdot (H_1+\epsilon H'_1) \cdot \ldots \cdot (H_{n-1}+\epsilon H'_{n-1}) < 0. $$
Then we find a sequence $(t_{\nu})$ of positive numbers converging to $0$ such that $T_X$ is stable w.r.t 
$(H_1+t_{\nu} H'_1,\ldots, H_{n-1}+t_{\nu} H'_{n-1})$ 
(adapt [CP07,5.5], where the case that all $H_i$ are equal is treated; actually by a convexity argument we have stability for all small $t$). 
Now we apply (1). 
\end{proof}  

For the next result compare also [RC00]. 

\begin{theorem} Let $X$ be a Fano manifold with $b_2(X) = 1.$ Then $T_X$ is generically ample (w.r.t. to $-K_X)$. 
\end{theorem} 

\begin{proof} Let $C \subset X$ be MR-general w.r.t. $-K_X$. Suppose that $T_X \vert C$ is not ample, hence
$T_X \vert C$ has a quotient of non-positive degree. Since $-K_X \cdot C > 0,$ we conclude that $T_X \vert C $ is not semi-stable. 
Consider the Harder-Narasimhan filtration of $T_X \vert C:$
$$ 0 = \sE_0 \subset \sE_1 \subset \ldots \sE_r = T_X \vert C.$$
By [KST07], there exists a positive $k$ such that $\sE_k$ is ample, and all ample subsheaves $\sF_C \subset T_X \vert C$ are contained in $\sE_k.$ In other words, $\sE_k$ is {\it maximally ample}. 
Consider the exact sequence
$$ 0 \to \sE_k \to T_X \vert C \to Q \to 0.$$
We claim that $c_1(Q) \leq 0. $ 
In fact, if $c_1(Q) > 0,$ then, being not ample, $Q$ must contain an ample subsheaf which contradicts the maximality of $\sE_k.$ 
Hence $$ c_1(\sE_k) \geq -K_X \cdot C. \eqno (*) $$
Now by [MR82], $\sE_k$ extends to a torsion free sheaf $\sE' \subset T_X$ of rank $k$, and by (*) 
$$ \det \sE' = -K_X \otimes \sO_X(a) \subset \bigwedge^k T_X $$
with an integer $a \geq 0$ (here $\sO_X(1)$ is the ample generator of ${\rm Pic}(X)$). 
Thus
$$ \sO_X(a) \subset \bigwedge^k T_X \otimes K_X = \Omega^{n-k}_X, $$ 
which is clearly impossible, e.g. by restricting to curves on which $T_X$ is ample (these exist by the rational connectedness of $X$). 
\end{proof}  

The proof of Theorem 5.4 actually shows

\begin{corollary}  Let $X$ be a projective manifold. \\
(1) Assume that $T_X$ is not generically ample w.r.t. $(H_1, \ldots, H_{n-1})$ and 
that $-K_X \cdot H_1 \cdot \ldots
\cdot H_{n-1} > 0.$ Then there exists a reflexive subsheaf $\sE \subset T_X$ such that 
$$c_1(\sE) \cdot H_1 \cdot \ldots \cdot H_{n-1} \geq -K_X \cdot  H_1 \cdot \ldots \cdot H_{n-1}$$
and such that $\sE \vert C$ is the maximal ample subsheaf, where $C$ is MR-general w.r.t. $(H_1, \ldots, H_{n-1}).$ \\
(2) Assume that  $T_X$ is not generically nef w.r.t. $(H_1, \ldots, H_{n-1})$ and 
that $-K_X \cdot H_1 \cdot \ldots
\cdot H_{n-1} > 0.$ Then there exists a reflexive subsheaf $\sE \subset T_X$ such that 
$$c_1(\sE) \cdot H_1 \cdot \ldots \cdot H_{n-1} > -K_X \cdot  H_1 \cdot \ldots \cdot H_{n-1}$$
and such that $\sE \vert C$ is the maximal nef subsheaf, where $C$ is MR-general w.r.t. $(H_1, \ldots, H_{n-1}).$ 
\end{corollary} 

\begin{proof} (1) is exactly what the proof of (5.4) gives. For (2), we argue with nef instead of ample bundles, using the
obvious nef version of [KST07,Prop.29].
\end{proof} 

\begin{remark} The same argument as in (5.4) also works for a variety $X$ which is $\bQ-$Fano with $\rho(X) = 1,$ i.e. $X$ is a 
normal projective $\bQ-$factorial variety
with at most canonical singularities. In fact, most things go over verbatim (of course one needs to deal with Weil $\bQ-$Cartier divisors), only at the end
the following additions need to be made. The variety $X$ is rationally connected by [Zh06]; we have a Weil divisor $\sL \subset\tilde  \Omega_X^{n-k}$ which is ample
or trivial, and the sheaf $\tilde \Omega^{n-k}_X$ is the extension of the sheaf of $(n-k)-$forms on the regular part of $X.$ Then we pass to a desingularization 
and argue there as before.
\end{remark} 

In the following we investigate the existence of {\it some} polarization $H$ such that $T_X$ is $H-$semi-stable 
resp.generically nef (ample).
\begin{proposition} Let $X$ and $Y$ be projective manifolds, $\pi: X \to Y$ birational. Assume that $T_Y$ is $(H_1, \ldots, H_{n-1})$-(semi-)stable.
Then $T_X$ is semi-stable w.r.t. $(\pi^*(H_1), \ldots, \pi^*(H_{n-1})).$ 
\end{proposition} 

\begin{proof} We argue for stability only. Write $H'_j = \pi^*(H_j)$ for short and $H' = (H'_1, \ldots, H'_{n-1});$ similarly introduce $H.$  
Assume that $T_X$ is $H'-$unstable and let $\sS'$ be maximally destabilizing. Then
$$ \mu_{H'}(\sS') \geq \mu_{H'}(T_{X}). $$ 
Let $\sS := \pi_*(\sS') \subset \pi_*(T_{X}) \subset T_Y.$
Now trivial calculations show
$$ \mu_{H'}(T_X) = \mu_H(T_Y) $$
and $$ \mu_{H'}(\sS') = \mu_H(\sS).$$
Thus we obtain
$$ \mu_H(\sS) \geq \mu_H(T_Y)$$
contradicting the stability of $T_Y.$ 
\end{proof} 

\begin{corollary} Let $\pi: X \to Y$ be a birational map of projective manifolds $X$ and $Y$. Assume that $T_Y$ is $(H_1, \ldots, H_{n-1})-$stable.
Then there are ample divisors $H'_j$ on $X$ such that $T_X$ is $(H'_1, \ldots, H'_{n-1})-$stable. 
\end{corollary}

\begin{proof} Combine the last proposition and [CP07,5.5].
\end{proof} 

\begin{corollary} 
Let $Y$ be a projective manifold with $(H_1, \ldots, H_{n-1})-$stable tangent bundle. Assume that 
$-K_Y \cdot H_1 \cdot \ldots H_{n-1} > 0$
(e.g. $Y$ is a Fano manifold of dimension 3 with $b_2 = 1 $ or  rational-homogenenous manifold
with $b_2 = 1$). Let $\pi: X \to Y$ be birational. Then there exists $H'_1, \ldots, H'_{n-1}$ ample on $X$ such that $T_X $ is 
$(H'_1, \ldots, H'_{n-1})-$stable and such that $T_X$ is generically nef w.r.t. $(H'_1, \ldots, H'_{n-1})$. 
\end{corollary} 

\begin{proof} By [CP07,5.5] we can take $H'_j = \pi^*(H_j) + \epsilon A$ where $A$ is ample and $\epsilon > 0$ small.
In order to prove generic nefness, we need to verify that
$$-K_X \cdot H'_1 \cdot \ldots \cdot H'_{n-1} > 0.$$
From $$-K_Y \cdot H_1 \cdot \ldots \cdot H_{n-1} > 0, $$ 
we deduce 
$$ -K_X \cdot \pi^*(H_1) \cdot \ldots \cdot \pi^*(H_{n-1}) > 0,$$ hence 
$$ -K_X  \cdot H'_1 \cdot \ldots \cdot H'_{n-1} > 0$$  
for $\epsilon $ small enough. 
\end{proof} 

\begin{corollary} Let $X$ be a rational surface. Then there exists $H$ ample such that $T_X $ is $H-$generically ample. 
\end{corollary} 

\begin{proof} Choose a birational morphism $\pi: \hat X \to X$ such that $\hat X$ is a blow-up of $\bP_2.$ By the last corollary we find
$\hat H$ ample on $\hat X$ such that $T_{\hat X}$ is $\hat H$-generically ample. Arguing by induction, we may assume that $\pi$ is just one blow-up;
let $E$ be the exceptional curve. Write 
$$\hat H = \pi^*(H) + aE$$
with $a < 0$ and $H$ ample on $X$. Let $\hat C \subset \hat X$ be MR-general w.r.t $\hat H.$ Then $T_{\hat X} \vert \hat C$ is ample. Since
$T_{\hat X} \subset \pi^*(T_X),$ also $\pi^*(T_X) \vert \hat C$ is ample. Therefore $T_X$ is ample on the image curve $C = \pi(\hat C).$ 
Since ampleness is a Zariski-open property, $T_X \vert C'$ is ample for $C'$ MR-general w.r.t. $H.$ 
\end{proof} 

Notice that for a Hirzebruch surface $\bP(\sO \oplus \sO(-e))$ with $e \geq 3,$ the tangent bundle is never stable for any polarization $H;$ however 
it is generically ample for all $H.$ \\
Corollary 5.10 can be strengthened as follows.

\begin{theorem} Let $X$ be a del Pezzo surface or a rational surface with $-K_X$ pseudo-effective. Let $H$ be ample. Then $T_X$ is $H-$generically ample. 
\end{theorem} 

\begin{proof} Suppose $T_X$ is not $H-$generically ample. By (5.5) there exists a line bundle $\sL \subset T_X$ such that 
$$ \sL \cdot H \geq -K_X \cdot H $$
and $\sL \vert C \subset  T_X \vert C$ is maximal ample for $C$ MR-general w.r.t. $H.$ Moreover $\sL$ is a subbundle outside a finite set. 
Notice that our assumptions guarantee $-K_X \cdot H > 0.$ 
By [KST07], [BM01], the general leaf of the foliation $\sL$ is algebraic and its closure is rational. So the closures of the (general) leaves form an algebraic 
family and we find a (smooth) blow-up 
$$ p: \sC \to X$$
which is a generic $\bP_1-$bundle $q: \sC \to B \simeq \bP_1$ such that 
$$ \sL = p_*(T_{\sC/B})^{**} = p_*(-K_{\sC/B})^{**} \subset p_*(T_{\sC}) = T_X.$$ 
Now
$$ \sL \cdot H = (-K_{\sC/B}) \cdot p^*(H) = -K_{\sC} \cdot p^*(H) + q^*(K_B) \cdot p^*(H) = $$
$$ = p^*(-K_X) \cdot p^*(H) +  q^*(K_B) \cdot p^*(H).$$
Since $B \simeq \bP_1,$ we have $ q^*(K_B) \cdot p^*(H) < 0$ which leads to
$$ \sL \cdot H < -K_X \cdot H, $$
a contradiction.  
\end{proof}

\begin{remark} {\rm The same scheme of proof stills works for birationally ruled surfaces over an elliptic curve with $-K_X$ pseudo-effective: \\
{\it the tangent bundle
of a surface with pseudo-effective anti-canonical bundle, which is birationally equivalent to a ruled surface over an elliptic curve, is generically nef.} \\
In the proof of (5.11) we take a priori $\sL \subset T_X$ such that $\sL_C$ is maximally nef. Since $-K_X \cdot C > 0$,
it is clear that $\sL_C$ is also maximally ample. Then we proceed as before, having in mind that $\sL \cdot H > -K_X \cdot H.$}
\end{remark} 

Combining Theorem 5.11 and Remark 5.12 we obtain

\begin{corollary} Let $X$ be a smooth projective surface with $-K_X$ nef. Then $T_X$ is generically nef. 
\end{corollary} 

\section{The movable cone} 
\setcounter{lemma}{0}

In order to generalize Theorem 5.4 we use the following cone theorem [BCHM07].
\begin{theorem} Let $X$ be a projective manifold with $-K_X$ big and nef. 
\begin{enumerate}
\item  $\sME^{-}(X) = \sME(X)$ is rationally polyhedral:   $\sME^{-}(X) = \sum_{j=1}^r R_j,$ where the $R_j$ are the extremal rays in $\sME(X).$ 
\item if $R \subset \sME^{-}(X)$ is an extremal ray, then there is a sequence of contractions and flips $X \dasharrow X'$ and a
Mori fiber space $\phi: X' \to Y$ such that the pull-back of a sufficiently general curve contracted by $\phi $ is in $R.$  
\end{enumerate}
\end{theorem} 

\begin{proof} In case $-K_X$ is ample, this is [BCHM07,1.3.4]. If $-K_X$ is merely big and nef, write
$$ -K_X = A + D$$
where $A$ is an ample $\bQ-$divisor and $D$ is an effective $\bQ-$divisor. 
Consider the birational morphism $\psi: X \to Y$
defined by $\vert -mK_X \vert $ for large $m.$ Then $-D$ is $\psi-$ample; consequently $-(K_X + \epsilon D)$ is ample for small positive
$\epsilon.$ Moreover the pair
$(X,\epsilon D)$ is klt, so that we can apply [BCHM07,1.3.4] to $K_X + \epsilon D.$ 
\end{proof} 

One should expect that $R$ contains the class of an irreducible rational curve whose deformations cover $X$. These curves should come from
suitable rational curves in a general fiber $F'$ of $\phi.$ The technical difficulty lies in the fact that $F'$ is singular (Fano) and one needs to
establish rational curves which avoid a set of codimension 2 in this singular variety.

Of course one could ask for more:
\begin{question}  Let $X$ be a projective manifold. 
\begin{enumerate} 
\item Is
$\sME^{-}(X)$ always locally rationally polyhedral?
\item If $-K_X$ is nef, is there a substitute of (6.1)(2), with log-contractions and log-flips?
\item If $X$ is rationally connected, is an extremal ray of $\sME^{-}(X)$ represented by a covering family of rational curves?
\end{enumerate} 
\end{question} 

Observe  however that (6.1)(2) without nefness assumption will in general not be true, even if $-K_X$ is big. Indeed, consider a 
ruled surface $X = \bP(E) $ over a curve  $C$ of genus $g \geq 2.$ 
Assume $X$ has
a negative section $C_0$ with $C_0^2 = -e < 0$ and $e = 2g(C)-1$ so that $-K_X$ is big. Let $F$ be a ruling line. Then the boundary components of $\overline{ME}(X)$ are 
represented by $F$ and $C_0 + eF.$ Then (6.1)(2) obviously does not hold, since the contraction of $C_0$ is not a Mori contraction; moreover the ray determined
by $C_0 + eF$ does not contain a covering family of (generically irreducible) curves. \\
If $X$ is an irrational surface with $-K_X$ nef, then obviously (6.2)(1) and (2) hold; here we have only one extremal ray.

\begin{theorem} Let $X$ be a projective manifold with $-K_X$ big and nef. 
Then $T_X$ is generically ample.
\end{theorem} 

\begin{proof} 
If the claim was false, we find by (5.5) some ample $(H_1, \ldots, H_{n-1}),$ and a subsheaf $\sE \subset T_X$ such that 
$$ (K_X + \det \sE) \cdot H_1 \cdot \ldots \cdot H_{n-1} \geq 0.$$
Notice that $\det \sE \subset \bigwedge^k T_X$ with $k$ the rank of $\sE$, and therefore, as in (5.4),
$$ K_X + \det \sE \subset \bigwedge^k T_X \otimes K_X = \Omega^{n-k}_X. \eqno (*)$$ 
We want to show that by a change of the polarizations $H_j$, we may achieve 
$$ (K_X + \det \sE) \cdot  H_1 \cdot \ldots \cdot H_{n-1} > 0. \eqno (**)$$
So assume equality in (**). 
Since the interior of $\sC\sI(X)$ is open in $N_1(X),$ we find ample $\bQ-$divisors $H_j'$ sufficiently near to $H_j,$ such that 
$$(K_X + \det \sE) \cdot  H'_1 \cdot \ldots \cdot H'_{n-1} > 0 $$
unless $K_X + \det \sE \equiv 0.$ But this last alternative cannot happen: since $\pi_1(X) = 0,$ we would have $K_X + \det \sE \simeq \sO_X$ and therefore by (*)
$h^0(\Omega_X^{n-k}) \ne 0$, contradicting $H^{n-k}(\sO_X) = 0$ (Kodaira vanishing).
So (**) can always be achieved.\\
By (6.1) there is an extremal ray $R \subset \sME^{-}(X)$ such that 
$$ (K_X + \det \sE) \cdot R > 0. $$
Moreover we can run a suitable MMP $X \dasharrow X'$ and obtain a Mori fiber space $\phi: X' \to Y$. Let $\sE' \subset T_{X'}$ be the induced
reflexive subsheaf, so that 
$$K_{X'} + \det \sE' \subset \tilde \Omega_{X'}^{n-k} $$ and
$$ (K_{X'} + \det \sE') \cdot R' > 0, $$
where $R'$ is the extremal ray contracted by $\phi.$ Therefore $ (K_{X'} + \det \sE')$ is a $\phi-$ample $\bQ-$Cartier divisor. Now we restrict to
a general fiber $F'$ of $\phi$ and, by the same techniques as fully explained in part (I) of the proof of the next theorem and by applying Lemma 6.6, we 
obtain a contradiction. 
\vskip .2cm \noindent
\end{proof}

One might ask whether the above method also works if $-K_X$ is merely nef, granted that Questions 6.2(1) and (2) have positive answers.
Of course we only ask for generic nefness. 
We try to argue by induction as before and this time obtain $\sE \subset T_X$ such that 
$$(K_X + \det \sE) \cdot H_1 \cdot \ldots \cdot H_{n-1} > 0.$$
So (**) automatically holds. The problem is now that there is no way to conclude the existence of an extremal ray $R$ in the movable cone such that
$$ (K_X + \det \sE) \cdot R > 0.$$ 
If we had $R$, then things go through smoothly. It is clear that we arrived here at a subtle point: up to now all arguments also work if $-K_X$ is merely
pseudo-effective, but we know that under this weaker assumption $T_X$ is in general not generically nef. \\
What we actually can conclude in case $-K_X$ nef is the existence of an element $\gamma \in \overline{ME}(X)$ such that
$$ K_X \cdot \gamma = 0 $$
and $$ (K_X + \det \sE) \cdot \gamma > 0. $$

\vskip.2cm \noindent
We continue to discuss the case that $-K_X$ is not big (but still nef). Here $-K_X$ might no longer be semi-ample. Geometrically speaking, we consider the
nef reduction  
$$\phi: X \dasharrow W $$ of $-K_X$ and assume 
that $0 < \dim W < \dim X = n$ (this is to say that $K_X \not \equiv 0$ and there is a covering family $(C_t)$ of curves such that
$K_X \cdot C_t = 0$), and that $\nu(-K_X) = \dim W.$ Then $\kappa (X) = \nu(X) $ and by Kawamata [Ka92] a suitable multiple
$-mK_X$ is spanned and thus $\phi$ is holomorphic and given by $\vert -mK_X \vert.$ Then $-K_X$ is automatically hermitian semi-positive, but we do not use
this information. The proof of Theorem 6.4 can also be achieved without using Kawamata's result, just using $\nu(-K_X) = \dim W.$ The reason is that
the almost holomorphicity of $\phi$ suffices to make our arguments work. We shall not use the finer structure of manifold with $-K_X$ hermitian semi-positive
[DPS96], see also (7.11) because we are aiming also for the general nef case; moreover the arguments in the proof of (6.4) also work in a singular
setting (canonical singularities). 

\begin{theorem} Let $X$ be a projective manifold with $-K_X$ semi-ample. Then $T_X$ is generically nef.
\end{theorem} 

\begin{proof} {\bf (I)} If $K_X \equiv 0,$ then Miyaoka's theorem [Mi87] already gives the result; hence we shall assume $K_X \not \equiv 0.$\\
 Let $\phi: X \to W$ be the fibration given by a suitable multiple $-mK_X$. By (6.3) we may assume $-K_X$ not to be big, hence
$\dim W < \dim X.$ 
Let $F$ be a general fiber of $\phi$ and set $d = n-k = \dim W.$ 
Then we have an exact sequence
$$ 0 \to \sO_F^d \to \Omega^1_X \vert F \to \Omega^1_F \to 0.  \eqno (1)$$ 
Assuming $T_X$ not to be generically nef, we apply (5.5)(2) and obtain a subsheaf $\sE \subset T_X$ of rank $k$
such that 
$$ (K_X + \det \sE) \cdot H_1 \cdot \ldots \cdot H_{n-1} > 0$$  
for suitable ample divisors $H_j.$ \\
We restrict the inclusion $K_X + \det \sE \subset \Omega^{n-k}_X$ to $F$ and obtain an injective map
$$ K_F + \det \sE_F \subset \Omega^{n-k}_X \vert F. \eqno (2)$$ 
We apply $\bigwedge^{n-k}$ to (1) and obtain sequences 
$$ 0 \to \sF_j \to \sF_{j-1} \to \sF_{j-1}/\sF_j \to 0 $$
with $$\sF_{j-1}/\sF_j \simeq \Omega^{n-k-j+1}_F \otimes \bigwedge^{j-1}\sO_F^d. $$
Now we chase the inclusion (2) through these sequences, using the fact that any morphism
$$ \psi: K_F + \det \sE_F \to (\Omega^r_F)^{\oplus N} $$
must vanish for $r > 0$. This will be proved below. The outcome is a non-vanishing morphism
$$ \psi: K_F + \det \sE_F \to \bigwedge^{n-k}\sO_F^d = \sO_F^N.$$ 
Such a morphism however cannot exist, since 
$$(K_X + \det \sE) \cdot B > 0$$ for a suitable movable curve $B \subset F,$ the existence of which also being proved
below. 
\vskip .2cm \noindent {\bf (II)} It remains to prove the vanishing of $\psi$ and the existence of $B.$ 
By projection to a suitable direct summand we obtain a non-zero morphism 
$$\tau: K_F + \det \sE_F \to \Omega^j_F, \eqno (3) $$
for some $j  \geq 0.$ 
Since $F$ is a fiber of $\phi,$ the morphism associated with $\vert -mK_X \vert,$ 
we have $K_F \equiv 0.$  
\vskip .2cm \noindent {\bf (III)} We want to reduce ourselves to the case that
$$ (K_X + \det \sE)_{> 0} \cap \overline {ME}(X) \cap K_X^{\perp} \ne \{0\}. \eqno (4) $$
Suppose (4) were false and let $H = (K_X + \det \sE)^{\perp}, $ so that 
$$\overline{ME}(X) \cap K_X^{\perp} \subset H_{\leq  0}. \eqno (4')$$ 
{\bf (III.1)} We choose a pseudo-effective $\bR-$divisor $P$ which is not big such that   
$$ P^{\perp} \cap \partial \overline{ME}(X) \subset H_{> 0} \eqno (5)$$ 
and such that $P$ is $\phi-$ample, i.e. positive on $K_X^{\perp}.$   
The divisor $P$ is constructed as linear form
$$\lambda: N_1(X) \to \bR $$ which is non-negative on $\overline{ME}(X)$ (using [BDPP04] to conclude that $P$ is 
pseudo-effective). More specifically, we apply Lemma 6.5 below with $ V = N_1(X); K = \overline{ME}(X) + \overline{NE}(X/W),$ the closed cone
generated by $\overline{ME}(X)$ and $\overline{NE}(X/W),$ 
with $L_1$ a linear subspace of $N_1(X) $ such that $L_1 \cap K \subset \partial K \cap H_{> 0}$ and $L_2$ the vector space generated by 
$\overline{NE}(X/W).$ The choice of $L_1$ is possible since $K_X + \det \sE$ is positive on parts of $\overline{ME}(X).$ \\
There is a slight difficulty with the application of Lemma 6.5, since we do not know that $L_2 \cap K = \overline{NE}(X/W) \subset H_{\leq 0}.$
Therefore we change $H$ to $\tilde H$ such that $K \cap H_{\leq 0} \subset K \cap \tilde H_{\leq 0}$ and such that $L_2 \cap K \subset \tilde H_{\leq 0}.$ 
Clearly we still can achieve $ \overline {ME}(X) \not \subset \tilde H_{\leq 0},$ since $\overline{ME}(X) \not \subset H_{\leq 0}.$ \\
Now we apply, working with $\tilde H,$ Lemma (6.5) to obtain a linear form $\lambda: N_1(X) \to \bR$ with the properties stated in (6.5). The form $\lambda$
is given by intersection with an $\bR-$divisor $P$ which is pseudo-effective by (6.5(1)) and [BDPP04]. Moreover $P$ is $\phi-$ample by (6.5)(2) and 
$P$ is not big by 6.5(3). (If we work with the original $H$, we would be forced to take $L_2 $ to the vector space generated by 
$\overline{ME}(X) \cap \overline{NE}(X/W)$ and could only conclude the relative bigness of $P$). Since $K \cap \tilde H_{> 0} \subset K \cap H_{> 0},$ 
claim (5) follows from 6.5(2).  
\vskip .2cm \noindent
{\bf (III.2)}    
Since $-K_X$ is the pull-back of an ample $\bQ-$divisor on $W$ and since $P$ is $\phi-$ample
we conclude that $A_0 = \lambda P - K_X$ is an ample $\bR-$divisor for $0 < \lambda \ll 1$.
Fix such a $\lambda$ so that $K_X + A_0 = \lambda P.$ Choose ample $\bQ-$divisors $A_k$ whose classes converge to the class of 
$\lambda P.$ Let $t_k$ be the effective threshold or pseudo-effective value, i.e., the smallest positive number such that 
$K_X + t_kA_k$ is pseudo-effective. Since $P$ is pseudo-effective but not big, the sequence $(t_k)$ converges to $1$ and therefore
$K_X + t_k A_k$ converges to $\lambda P$ (numerically). From (5) we deduce 
$$ (K_X + t_0A)^{\perp} \cap \partial \overline{ME}(X) \subset H_{> 0} \eqno (6)$$
for $k \gg 0.$ We fix such a large $k$, choose $m$ such that $A := mA_k$ is Cartier and set $t_0 = mt_k.$ 
By [BCHM06,1.1.7] $t_k \in \bQ_+$ and so does $t_0.$ \\
For $\epsilon > 0$ rational and sufficiently small the $\bQ-$divisior $K_X + (t_0 - \epsilon)A $ is not pseudo-effective, 
so by [BCHM06] we can run the MMP and obtain a birational rational
map 
$$\sigma: X \dasharrow X'$$ 
(composed by birational contractions and flips) together with a Mori contraction 
$$f: X' \to Y $$ 
such that $\dim Y < \dim X.$ Notice that  
there is an algebraic set $B \subset X'$ of codimension at least 2 such that $\sigma$ is an isomorphism over $X' \setminus B.$ 
Let $A'$ be the divisor on $X'$ induced by $A.$ Then by construction $K_{X'} + t_0A'$ is the pull-back of an ample $\bQ-$divisor on $Y.$ 
Let $M'$ be the divisor on $X'$ induced by $\det \sE$ and let $l'$ be a general curve in $F',$ a general fiber of $f.$ Then $l' \cap B = \emptyset$ and
we can consider its isomorphic preimage $l \subset X.$ We claim that
$$ (K_{X'} + M') \cdot l' > 0. \eqno (7) $$ 
To prove (7), we need to verify $(K_X + \det \sE) \cdot l > 0,$ i.e. $[l] \in H_{> 0}.$ By (6) this holds as soon as we know that
$$ (K_X + t_0A) \cdot l = 0,$$
which however is obvious via 
$$(K_X + t_0A) \cdot l = (K_{X'} + t_0A') \cdot l' = 0.$$ 
So (7) holds and thus $K_{X'} + M' \vert F'$ is ample (since $\rho(X'/Y) = 1). $
Now, taking $\bigwedge^{n-k}$ of the exact sequence of K\"ahler differentials on the regular locus of the normal variety $F'$ 
$$ N^*_{F'\vert X'} \to \Omega^1_{X'} \vert F' \to \Omega^1_{F'} \to 0,$$
we conclude via Lemma 6.6 that $K_{X'} + M' \vert F'$ is a subsheaf of some $\tilde \Omega^{j}_{F'}. $ 
This is impossible, again by Lemma 6.6.
\vskip .2cm \noindent 
{\bf (IV)} 
So (4) holds and consequently there exists $\gamma \in \overline{ME}(X) \setminus \{0\}$ with $K_X \cdot \gamma = 0$ such that
$$ (K_X + \det \sE) \cdot \gamma > 0.$$ 
By (6.7) there exists $\delta \in \overline {ME}(F)$ such that 
$$ \gamma = i_*(\delta).$$ 
By pertubation, we may assume $\delta$ to be in the interior of $ME(F).$ Moreover we have a movable curve 
$$ B = \mu_*(H'_1 \cdot \ldots \cdot H'_{n-2}) $$
(for some modification $\mu: \tilde F \to F)$ such that 
$$(K_F + \det \sE_F) \cdot B > 0. $$ 
Since $K_F \equiv 0$ and since $K_F + \det \sE_F \subset \Omega^j_F$ by (3), this contradicts [CP07] (and is actually trivial for $j = 0$). 
\end{proof}

We still need to establish the following lemmata 6.5 and 6.6 and theorem 6.7.
\begin{lemma} Let $V$ be a finite-dimensional real vector space, $K \subset V$ a closed cone with $K \cap -K = \{0\}$ and $L_1,L_2 \subset V$ 
linear subspaces such that $L_1 \cap K \subset \partial K.$ Let $H \in V^* \setminus \{0\}$ such that $L_2 \cap K \subset H_{\leq 0}$ and $L_1 \cap K \subset H_{> 0}.$
Then there exists a non-trivial linear form $\lambda: V \to \bR$ such that
\begin{enumerate}
\item $\lambda \vert K \geq 0 $
\item $\lambda \vert K \cap H_{\leq 0} > 0$
\item $\lambda \vert L_2 \cap K > 0$    and 
\item $\lambda \vert L_1 = 0.$
\end{enumerate} 
\end{lemma} 

\begin{proof} Let $K_i$ be a sequence of closed cones with the following properties
\begin{itemize}
\item  $K_i \subset K_{i+1},$
\item $K_i \cap L_1 = \{0\} $ for all $i$ 
\item $K = \bigcup K_i$ and
\item $K_i \cap H_{\leq 0} = K \cap H_{\leq 0} $ for all $i$.
\end{itemize} 
By Hahn-Banach we may choose $\lambda_i: V \to \bR$ linear such that $\lambda_i \vert K_i > 0,$ $\lambda_i \vert L_1 = 0$ and $\Vert \lambda_i \Vert = 1
$ (for some norm in $V^*$). Then we extract a limit $\lambda: V \to \bR$, and automatically (1) and (4) hold. 
As to (2), a priori $\lambda $ could have zeroes on $K \cap H_{\leq 0}. $ This can be avoided by working from the beginning with some suitable $\tilde H$ 
sufficiently near to $H$ such that $K \cap H_{\leq 0}$ is in the interior of $K \cap \tilde H_{\leq 0}.$ \\
(3) is finally a consequence of (2).   
\end{proof} 
\begin{lemma} Let $X$ be a $\bQ-$Fano variety, i.e. a normal projective $\bQ-$Gorenstein variety with at most terminal singularities such that
$-K_X$ ample. Let $\sL \subset \tilde \Omega^r_X$ be a reflexive rank 1 subsheaf which $\bQ-$Cartier. Then $\kappa (\sL) = - \infty$.
\end{lemma}

\begin{proof} Suppose $\kappa (\sL) \geq 0$ and take an integer $m$ such that $h^0(\sL^{\otimes m}) \ne 0.$ Then
$$ H^0(X,((\tilde \Omega^r_X)^{\otimes m})^{**}) \ne 0.$$
Thus we find a positive integer $N$ such that
$$ H^0(X,((\tilde \Omega^1_X)^{\otimes N})^{**}) \ne 0.$$  Let $\pi: \hat X \to X$ be a desingularization.
Since $\pi_*(\Omega^1_{\hat X}) = \tilde \Omega^1_X$ by [Fl88], we obtain a generically injective map
$$ \pi_*[(\tilde (\Omega^1_X)^{\otimes N})^{**}] \to (\Omega^1_{\hat X})^{\otimes N},$$
and therefore $$ H^0(\hat X,(\Omega^1_{\hat X})^{\otimes N}) \ne 0. $$
But $\hat X$ is rationally connected by [Zh06], which is absurd (e.g. by [Ko95,IV.3.8]). 
\end{proof} 
For future use we state the following Theorem 6.7 in more general form than necessary for our purposes here. 
\begin{theorem} 
Let $X$ be a projective manifold, $L$ a nef line bundle on $X$ with (almost holomorphic) nef reduction $f: X \dasharrow W.$ 
Let $F$ be a general fiber of $f$ and $i: F \to X$ the inclusion map. We assume that 
\begin{enumerate} 
\item $L \vert F = \sO_F$;
\item $\dim W = \nu(L). $
\end{enumerate} 
Let $\gamma \in \overline{ME}(X) $ with $L \cdot \gamma = 0.$ Then there exists $\delta \in \overline {ME}(F)$ such that 
$\gamma = i_*(\delta).$
\end{theorem}

\begin{proof} 
(1) In a first step we show that we may assume $f$ holomorphic. In fact, take a birational map $\pi: \hat X \to X$ from a projective manifold
$\hat X$ such that the induced map $\hat f: \hat X \to W$ is holomorphic. Of course, we blow-up only inside the set of indeterminacies of $f.$
Then $\hat f$ is a nef reduction for $\hat L = \pi^*(L)$ and obviously $\nu(\hat L) = \nu(L).$ Next observe that 
$$ \hat \gamma := \pi^*(\gamma) \in \overline{ME}(\hat X). $$
To see that, take a pseudo-effective divisor $\hat D$. Then the divisor $\pi_*(\hat D)$ is pseudo-effective, too,  and therefore
$$\hat D \cdot \hat \gamma = \pi_*(D) \cdot \gamma \geq 0.$$
Now apply [BDPP04] to conclude. \\
Since $\hat L \cdot \hat \gamma = 0$ and since we assume in this step the assertion of the proposition to be true in case of a holomorphic
nef reduction, we find an element $\hat \delta \in \overline{ME}(\hat F)$ (where $\hat F$ is a general fiber of $\hat f$) such that 
$$ i_*(\hat \delta) = \hat \gamma. $$
But $F = \hat F,$ since $f$ is almost holomorphic, and thus our claim follows.
\vskip .2cm \noindent (2) From now on we assume $f$ holomorphic. By a birational base change and the above arguments, we may also 
assume $W$ smooth and $f$ flat, in particular equidimensional.  However $X$ is now only normal. 
\vskip .2cm \noindent
Consider the line bundle
$$ A = f_*(L)^{**},$$
here we use the assumption $L \vert F = \sO_F$ to see that $A$ has indeed rank $1$.
The canonical map $f^*f_*(L) \to L $ gives rise to a map $\psi: f^*(A) \to L$ which is defined outside of the preimage of the singular locus 
$S$ of the torsion free sheaf $f_*(L).$ But $S$ has codimension at least $2$, so does $f^{-1}(S),$ since $f$ is equidimensional. Therefore $\psi$ is defined everyhwere,
and we can write
$$ L = f^*(A) + E$$
with an effective divisor $E$. Since $E$ is $f-$nef with $f_*(\sO_X(E)) = \sO_W$, it is easy to show [Fu86,1.5] that 
$$ mE = f^*(B) $$
for some multiple $m.$ In total we can write
$$ mL = f^*(A') $$
with some nef line bundle $A'$. Notice that $A'$ also big by our assumption $\dim W = \nu(L).$
By passing to a desingularization of $X,$ we may again assume that $X$ is smooth; the new map $X \to W$ might no longer be equidimensional,
but this is not important for the rest of the considerations. 
\vskip .2cm \noindent   (3)  We show that 
$$ f_*(\gamma) = 0.$$
By (2),
$$ 0 = L \cdot \gamma = A' \cdot f_*(\gamma). $$ 
Since $f_*(\gamma) \in \overline{ME}(W) $ (again test by intersecting with a pseudo-effective divisor on $W$), 
we conclude from the bigness of $A'$ that $f_*(\gamma) = 0.$ 
\vskip .2cm \noindent (4) Now that we know $f_*(\gamma) = 0$, we apply Theorem 6.8 to conclude.
\end{proof}

\begin{theorem} Let $f: X \to W$ be a surjective holomorphic map with connected fibers from the projective manifold $X$ to the 
normal projective variety $W.$ Let $F$ be a general fiber of $f.$ Then 
\begin{enumerate} 
\item The inclusion map $i_*: N_1(F) \to N_1(X)$ maps $\overline {ME}(F)$ to $\overline {ME}(X).$ 
\item $\overline{ME}(X/W)$ is the closed cone generated by covering families $(C_t)$ in $X$ such that $f_*(C_t) = 0.$ 
\item $i_*: \overline {ME}(F) \to \overline {ME}(X/W)$ is surjective. 
\end{enumerate} 
\end{theorem}  

\begin{proof} (1) By [BDPP04], $\overline{ME}(F)$ is the closed cone generated by (the classes of the) irreducible curves $C$ with nef normal sheaves $N_{C/F}$. 
Here the normal sheaf $N_{C/F}$ is defined as
$$N_{C/F} = (\sI_C/\sI_C^2)^*,$$ 
where $\sI_C$ is the defining ideal. Recall also that a coherent sheaf $\sS$ is nef by definition, if the line bundle 
$\sO_{\bP(\sS)}(1)$ is nef on $\bP(\sS).$ 
Therefore we only need to show that 
$$i_*(C) \in \overline {ME}(X)$$ holds for a curve
$C$ with $N_{C/F}$ nef. To do that we let $\sI'$ be the ideal of $C \subset X$ and $\sJ$ be the ideal of $F \subset X.$ Then we have an exact sequence
$$ 0 \to (\sJ/\sJ^2) \vert C \to \sI'/\sI'^2 \to \sI/\sI^2 \to 0.$$
Now the morphism  $(\sJ/\sJ^2) \vert C = \sO_C^d \to \sI'/\sI'^2$ clearly does not have zeroes, hence we dualize to obtain an exact sequence
$$ 0 \to N_{C/F} \to N_{C/X} \to \sO_C^d \to 0.$$
Thus $N_{C/X}$ is nef and $C$ is a movable curve in $X.$ 
\vskip .2cm \noindent
(2) Let $L$ be a line bundle such that 
$$L \cdot C_t  \geq 0$$
for all covering families $(C_t)$ in $X$ with $f_*(C_t) = 0$. We need to show that
$$ L \cdot \gamma \geq 0 \eqno (*) $$
for  $\gamma \in \overline {ME}(X/W).$  \\ 
By (1), $L_F$ is pseudo-effective for the general fiber $F.$ Choose an ample line bundle $A$ on $X$. 
Then $L + \epsilon A$ is $f-$big for all $\epsilon > 0.$ Hence
$$ (L + \epsilon A) \cdot \gamma > 0,$$
and (*) follows by letting $\epsilon $ converge to $0.$
\vskip .2cm \noindent 
(3) is a consequence of (2). 
\end{proof} 

\begin{remark} {\rm  Suppose $-K_X \cdot C > 0$ for all movable curves $C$. One might expect that then $q(X) = 0$ or even $X$ to be rationally connected. 
One reason for such an expectation is Serrano's theorem, that if $-K_X \cdot C > 0$ for {\it all} curves in a threefold $X$, then $X$ is Fano (the same
being expected in all dimensions). This expectation however is completely false, even if $-K_X$ is nef. 
Take an elliptic curve $C$ and a rank 2-vector bundle $E$ given by a non-split
extension 
$$ 0 \to \sO_C \to E \to \sO_C \to 0. $$
Let $X = \bP(E)$. Then $-K_X \cdot B = 0$ only for one curve, the section defined by the epimorphism $E \to \sO_C \to 0.$ }
\end{remark} 

We next prove a weak substitute for generic nefness for rationally connected manifolds. To do that we first observe

\begin{theorem} Let $X$ be a rationally connected projective manifold. Then there are covering families $(C_t^j)$ of rational curves, $1 \leq j \leq m$,
such that $\sum_{j} [C_t^j]$ is in the interior ${ME}^0(X).$
\end{theorem}  

\begin{proof} 
Applying [Ko95,IV.3.9] we find for all $x_1,x_2\in X$ an irreducible rational curve $C_{x_1,x_2}$ joining $x_1$ and
$x_2$. By the usual Chow scheme argument, there are finitely many connecting families  $(C_t^j)$ of rational curves, 
with $T_X \vert C_t^j$ ample for general $t$ and all $j$,
$1 \leq j \leq m$, such that any two $x_1, x_2$ can be joined an irreducible member of some of the families. Moreover at least one of the families, say $(C_t^1)$,
has the property that two general points of $X$ can be joined by an irreducible member of this family. We may also assume that 
for any $j$ the general $C_t^j$ is smooth (resp. has only nodes when $X$ is a surface). 
In order to prove that $\alpha := \sum_{j} [C_t^j]$ is in the interior ${ME}^0(X),$ we check that
$$ L \cdot \alpha > 0$$
for all pseudo-effective line bundles $L$ on $X$ which are not numerically trivial. Suppose to the contrary that
$ L \cdot \alpha = 0$ for some $L$,
so that 
$$L \cdot C^j_t  = 0 \eqno (*) $$ for all $j$, in particular
$$ L \cdot C^1_t = 0.$$
Let $L = P + N$ be the divisorial Zariski decomposition [Bo04] of $L$; see also [BDDP04] for a discussion. So $P$ is an effective $\bR-$divisor
and $N$ is nef in codimension $1.$ Thus
$$P \cdot C_t^1 = N \cdot C_t^1 = 0 $$ and by [BDPP04,8.7], we first obtain $N = 0$ and consequently [ibid.,3.11] that $P$ is numerically equivalent 
to a multiple of an effective divisor, hence w.l.o.g. effective.
Let $P_k$ be the irreducible components of the support of $P.$ Since $P \cdot C^j_t = 0$ and since 
$P_k \cdot C^j_t \geq 0$ for all $j$, we have $P_k \cdot C^j_t = 0$ for all $j.$ Now pick $x_1$ general and $x_2 \in P_k$ and choose some
irreducible $C_t^j$ joining $x_1$ and $x_2.$ Thus $0 < P \cdot C_t^j = L \cdot C_t^j$ contradicting (*). 
\end{proof}

It seems likely that we can find a single ``connecting'' family $(C_t)$ of rational curves (with $C_t$ irreducible!) such that $[C_t]$ is in the 
interior of the movable cone, but we can't prove that at the moment. \\ 
As a consequence of Theorem 6.10 we obtain
\begin{theorem} Let $X$ be a rationally connected manifold. Then there exists a smooth curve $C \subset X$ such that
\begin{enumerate}
\item the deformations of $C$ cover $X;$
\item $T_X \vert C$ is ample;
\item $[C] \in ME^{o}(X).$
\end{enumerate} 
\end{theorem} 

\begin{proof} Choose families $(C_t^j)$ on $X$ as in Theorem 6.10. We form a reducible connected curve
$$ C' = \sum C_t^j, $$
where the $C_t^j$ are choosen general so that the intersections are all transversal and that through any point of $X$ there are at most two components.
Then $C'$ is smoothable ([Ko95,I.6.10]); let $C$ be a smoothing so that $[C] \in (ME(X))^{o}.$ 
Since $T_X \vert C_t^j$ is ample for all $j$ and $t$ general, so does $T_X \vert C$, proving (2). 
Claim (1) is also clear since $[C']$ belongs to a component of the Chow scheme which covers $X.$ \\


\end{proof}

The fact that  $[C] \in ME^{o}(X)$ can be rephrased as follows. There are birational morphisms $\pi_j: \tilde X_j \to X$ from projective manifolds
$\tilde X_j$ and $\tilde h_j = (\tilde H_1^j, \ldots, \tilde H^j_{n-1})$ ample such that $[C] = \sum_j a_j \pi_{j*}(\tilde h_j)$ with positive 
numbers $a_j.$ It is not clear whether one can simply acieve $[C] = \pi_*(\tilde H \cdot \ldots \cdot \tilde H_{n-1}). $\\

\section{Base Points} 
\setcounter{lemma}{0}

Generic nefness is sometimes not good enough for applications, e.g. for the smoothness of the Albanese map. Therefore we need to consider nefness concepts
which deal with curves through any point.

\begin{theorem} Let $E$ be a vector bundle on the projective manifold $X_n$ and $x \in X.$ Assume that $E$ is $(H_1, \ldots, H_{n-1})-$stable. 
Then there exists a covering family $(C_t)$ through $x$ such that $E \vert C_t$ is stable for
general $t$, i.e. $f_t^*(E)$ is stable, where $f_t$ is the normalization of $C_t.$  
\end{theorem} 

\begin{proof} Let $\pi: \tilde X \to X$ be the blow-up of $x$ with exceptional divisor $D$ and fix an ample divisor $A$ on $\tilde X.$ 
Then for any $\epsilon > 0,$ the divisor
$$ \tilde H_i = \pi^*(H_i) + \epsilon A $$
is ample for all $i$ 
and $\pi^*(E)$ is stable w.r.t. $(\pi^*(H_1), \ldots, \pi^*(H_{n-1})).$ \\
Hence by [CP07,5.5] (with a trivial modification) $\pi^*(E)$ is $(\tilde H_1, \ldots, \tilde H_{n-1})-$stable for suitable sufficientl small $\epsilon$.
Now we apply Mehta-Ramanathan:
let $\tilde C$ be MR-general w.r.t.  $(\tilde H_1, \ldots, \tilde H_{n-1}),$ so that $\pi^*(E) \vert \tilde C$ is stable. Since $\pi \vert \tilde C $ is
the normalization of $\pi(\tilde C)$ and since $D \cap \tilde C \ne \emptyset,$ so that $x \in \pi(\tilde C)$, the proof is complete.
\end{proof} 

\begin{corollary} Suppose the vector bundle $E$ to be  $(H_1, \ldots, H_{n-1})-$stable and that $\det E \cdot H_1 \cdot \ldots \cdot H_{n-1} \geq 0.$ 
Then $E$ is sufficiently nef. If $\det E \cdot H_1 \cdot \ldots \cdot H_{n-1} > 0,$ then $E$ is sufficiently ample. 
\end{corollary} 

\begin{proof} Fix $x \in X.$ By (7.1) there exists a covering family $C_t$ through $x$
such that $E \vert C_t$  is stable for general $t.$ Since $\det E \cdot H_1 \cdot \ldots \cdot H_{n-1} \geq 0$, we conclude as in the proof
of (5.1) that $E \vert C_t$ is nef, analogously for ampleness.
\end{proof}

\begin{corollary} Let $X$ be a projective manifold with $K_X \equiv 0.$ Then $T_X$ is sufficiently nef. 
\end{corollary} 

\begin{proof} Fix ample divisors $H_i.$ A priori $T_X$ is only semi-stable w.r.t. any $(H_1, \ldots, H_{n-1}).$ Then there is a decomposition (see e.g. [Pe94])
$$ X \simeq \Pi X_j $$
such that $T_{X_j} $ is stable for all $j.$ Hence $T_{X_j}$ is sufficiently nef for all $j$ and so does $p_j(T_{X_j}).$ Since 
$$ T_X = \bigoplus p_j^*(T_{X_j}),$$
where $p_j: X \to X_j$ are the projections, we conclude easily that $T_X$ is sufficiently nef.  
\end{proof} 

\begin{remark} {\rm  The arguments of (7.3) actually show the following. If $X \simeq \Pi X_j$ with all $T_{X_j}$ sufficiently nef, then $T_X$ is
sufficiently nef.}
\end{remark} 

Corollary 7.2 implies that the tangent bundle of ``many'' Fano manifolds is sufficiently nef, even sufficiently ample. In the following we improve 
(7.2) and obtain in particular  that the tangent bundle of {\it any}
Fano manifold is sufficiently ample - this could also be deduced from rational connectedness (at least for curves through general $x$).  


\begin{theorem} Assume that the vector bundle $E$ is generically ample w.r.t. \\ $(H_1, \ldots, H_{n-1}) $. 
Let $x \in X$. Then there exists a covering family $(C_t)$ through $x$ such that $E \vert C_t$ is ample for general $t$, i.e., $E$ is
sufficiently ample.
\end{theorem} 

\begin{proof} Let $\pi: \tilde X \to X$ be the blow-up of $x$. Fix an ample divisor $A$ on $\tilde X$ and 
write
$$ \tilde H_{i,\epsilon} = \pi^*(H_i) + \epsilon A. $$  
For $\epsilon > 0$ rational and sufficiently small, $ \tilde H_{i,\epsilon}$ is an ample $\bQ-$divisor. 
Assuming our claim to be false, we conclude that for $\tilde C \subset \tilde X$ MR-general w.r.t. $ \tilde H_{1,\epsilon}, \ldots,  \tilde H_{n-1,\epsilon}$
the bundle
$$ \pi^*(E) \vert \tilde C $$
is not ample. 
Let $$\sS_{\epsilon} \subset \pi^*(E) \vert \tilde C$$
be maximally ample.
Then as in (5.4) $\sS_{\epsilon}$ extends to a global subsheaf 
$$\tilde \sS_{\epsilon} \subset \pi^*(E).$$
Notice that $\tilde \sS_{\epsilon}$ is a part of the Harder-Narasimhan filtration w.r.t  $ \tilde H_{1,\epsilon}, \ldots,  \tilde H_{n-1,\epsilon}$, moreover
by maximal ampleness we have
$$ c_1(\tilde \sS_{\epsilon}) \cdot  \tilde H_{1,\epsilon} \cdot \ldots \cdot  \tilde H_{n-1,\epsilon}  \geq
\pi^*(\det E)  \cdot  \tilde H_{1,\epsilon} \cdot \ldots \cdot  \tilde H_{n-1,\epsilon}. \eqno (*) $$
By (the proof of) [CP07,5.5], the sheaf $\tilde \sS_{\epsilon}$ does not depend on $\epsilon$ for $\epsilon $ suitable and sufficiently small; 
we therefore drop the index (the proof of [CP07,5.5] actually gives that for the maximal destabilising subsheaf; otherwise consider the quotient
by the maximal destabilizing subsheaf and proceed by induction).  
Now we pass to the limit $\epsilon \to 0$ in (*) and
obtain
$$ c_1(\tilde \sS) \cdot \pi^*(H_1) \cdot \ldots \cdot \pi^*(H_{n-1}) \geq \pi^*(\det E)  \cdot \pi^*(H_1) \cdot \ldots \cdot \pi^*(H_{n-1}).$$
Let
$$\sS = \pi_*(\sS) \subset E.$$ Then 

$$ (\det E^* + \det\sS) \cdot H_1 \cdot \ldots \cdot H_{n-1} \geq 0. $$
On the other hand, $\det E^* + \det \sS \subset \bigwedge^{p}E^*$  for a suitable $p$ and thus we obtain a contradiction to the assumption that 
$E$ is generically ample w.r.t. $(H_1, \ldots, H_{n-1}).$

\end{proof} 

By (6.3) we conclude

\begin{corollary} If $-K_X$ is big and nef, then $T_X$ is sufficiently ample. 
\end{corollary} 

\begin{remark} {\rm The analogous statement: 
\vskip .1cm \noindent 
\centerline{{\it  $E$ generically nef implies $E$ sufficiently nef }}
\vskip .1cm \noindent 
is unfortunately false. E.g. let $X = \bP_3$ and 
consider a rank 2-vector bundle $\sE$ on $\bP_3$ given as an extension 
$$ 0 \to \sO(a) \to \sE \to \sI_Z \to 0 $$
with a suitable $a > 0$ and a suitable locally complete intersection curve $Z \subset \bP_3$. Then obviously $\sE$ is generically nef, 
since the general curve does not meet $T$, but neither generically ample nor sufficiently nef ($\sE \vert C$ is never nef for a curve $C$ which meets
$Z$ is a finite set). \\
The difficulty with going to the limit in (7.5) is the following.  Fix $x \in X$ and an ample line bundle $\sL.$ We now apply (7.5) to the $h-$generically ample bundle $S^mE \otimes \sL.$ So if
$x \in X$ we find a covering family $(C_t)$ through $x$ such that $S^m \otimes \sL \vert C_t$ is ample for general $t.$ 
The family $(C_t)$ however depends on $m$ and therefore one cannot conclude. } 
\end{remark} 

\begin{remark} {Although (7.5) fails in the generically nef case, the proof of (7.5) still gives something. 
So suppose that $E$ is generically nef  w.r.t. $(H_1, \ldots, H_{n-1}) $, the same being also true for small pertubations of $H_j.$ 
Suppose furthermore that sufficient nefness fails at $x$. 
Following the lines of argumentation in (7.5) (substituting maximal ampleness by maximal nefness), 
we obtain
$$ \tilde \sS \subset E$$
such that
$$ c_1(\tilde \sS) \cdot \pi^*(H_1) \cdot \ldots \cdot \pi^*(H_{n-1}) \geq \pi^*(\det E)  \cdot \pi^*(H_1) \cdot \ldots \cdot \pi^*(H_{n-1}).$$
We again consider $\sS = \pi_*(\tilde \sS) \subset E$ so that
$$ (\det E^* + \det \sS) \cdot H_1 \cdot \ldots \cdot H_{n-1} \geq 0. \eqno (*) $$
If we have strict inequality in (*), then we can conclude as in (7.5) to get a contradiction. Thus 
equality in (*) must happen. Then for small pertubations of $H_j$, we also cannot have a strict inequality (*);
consequently 
$$ \det E^* + \det \sS \equiv 0.$$
Introducing $Q = T_X/\sS,$ we obtain $\det Q \equiv 0.$ Moreover we find a covering family $(C_t)$ through $x$ such that $Q \vert C$ is negative modulo torsion. 
To see this, consider $\tilde Q = \pi^*(E)/ \tilde \sS$ and observe that by construction 
$$ \tilde Q \vert \tilde C $$
is negative, where $\tilde C$ is MR-general  $ \tilde H_{1,\epsilon}, \ldots,  \tilde H_{n-1,\epsilon}$ for suitable small $\epsilon.$ 
We obtain an exact sequence 
$$ 0 \to \pi_*(\tilde \sS) \to T_X \to \pi_*(\tilde Q) \to R^1\pi_*(\tilde \sS) \to 0,$$
so that $Q$ is a subsheaf of $\pi_*(\tilde Q)$ which coincides with $\pi_*(\tilde Q)$ outside $x.$ Since $\pi^*\pi_*(\tilde Q)$ is a 
subsheaf of $\tilde Q$ modulo torsion, 
$$ \pi_*(Q) \vert \pi(\tilde C)$$
is negative and so does $Q \vert \pi(\tilde C).$}
\end{remark} 

\begin{proposition}
Let $X$ be a rationally connected manifold. If $T_X$ is generically nef, $T_X$ is sufficiently nef. 
\end{proposition} 

\begin{proof} We assume that $T_X$ is not sufficiently nef at $x.$ By (7.8) we obtain a sequence
$$ 0 \to \sS \to T_X \to Q \to 0$$
with $\det Q \equiv 0.$ Since $-K_X$ is big and nef, $\det Q = \sO_X,$ hence 
$$ H^0(X,\Omega^j_X) \ne 0$$
for some $j \geq 1$ a contradiction. 
\end{proof}

The same proof as in (7.9) actually shows the following (using the almost holomorphicity of the rational quotient):

\begin{proposition} Suppose $T_X$ generically nef, but not sufficiently nef. Let $f: X \dasharrow W$ be the rational quotient and $Q$ the quotient
sheaf of $T_X$ constructed in (7.8). Then the composed rational map $T_{X/W} \to Q$ vanishes; in particular ${\rm rk}Q \leq \dim W.$ 
\end{proposition}

\begin{theorem} Suppose that $-K_X$ is hermitian semi-positive. 
\begin{enumerate}
\item  The rational quotient has a model which is a holomorphic fiber bundle $f: X \to Z$. Possibly after finite \'etale cover, 
$K_Z = \sO_Z $ and $Z \simeq A \times B$ with $A$ the Albanese torus of $X$ and $B$ a product of hyperk\"ahler and Calabi-Yau manifolds. 
\item  The tangent bundle $T_X$ is sufficiently nef. 
\end{enumerate}
\end{theorem} 

\begin{proof}(1) By [DPS96], and [Zh05] for the rational connectedness statement, the Albanese map 
$$ \alpha: X \to A$$
is a fiber bundle over $A$ and possibly after a finite \'etale cover of $X,$ the fiber $F$ of $\alpha$ is a product 
$$ F = \Pi F_j$$ 
such that the $F_j$ are hyperk\"ahler, Calabi-Yau or rationally connected. From that it is immediately clear that the 
rational quotient has a model which is a holomorphic fiber bundle $X \to Z.$ Moreover $Z $ is holomorphic fiber bundle
over its Albanese such that the fibers $F$ have $K_F \equiv 0.$ Therefore $K_Z \equiv 0$ and after a finite \'etale cover
$Z$ is a product as stated.
\vskip .2cm \noindent 
(2) The tangent bundles $T_{F_j}$ are sufficiently nef by (7.3) resp. by (6.4) and (7.9). Therefore $T_F$ is sufficiently
nef by (7.4). This already solves the problem in case $\pi_1(X)$ finite.  
\vskip .2cm \noindent 
Now assume that $T_X$ is not sufficiently nef at $x$ and apply (7.8).
In the following we will be free to pass to a finite \'etale cover if necessary.
Consider the canonical morphism
$$ \lambda: T_{X/A} \to T_X \to Q$$
and set 
$$ Q' = {\rm Im} \lambda; \ \sS' = {\rm Ker} \lambda.$$
Notice that the map $\tilde \sS \to \pi^* \alpha^*(T_A) = \sO_{\tilde X}^d $ 
has not maximal rank at most in codimension 2  (restrict to MR-general curves and use nefness of $\tilde \sS$ on those curves), 
Now an easy diagram chase shows that $Q'$ has the same negativity property as $Q$ in (7.8) since there is an exact sequence
$$ 0 \to Q' \to Q \to {\mathcal T} \oplus \sO_X^e \to 0 $$
with ${\mathcal T}$ a torsion sheaf; supported in codimension at least 2 by the sentence before and some $e \geq 0.$ 
In particular
$$ \det Q' \equiv 0.$$
Now consider the relative tangent bundle sequence
$$ 0 \to T_{X/Z} \to T_{X/A} \to f^*(T_{Z/A}) = p_B^*(T_B) \to 0.$$ 
By restricting to a fiber of $f$ (which is rationally connected) and by recalling $\det Q' \equiv 0$ the
induced map 
$$ \lambda: T_{X/Z} \to Q' $$
vanishes. In fact we
consider the exact sequence
$$ 0 \to {\rm Im} \lambda_F \to Q' \to {\rm Coker} \lambda_F \to 0.$$
Since ${\rm Coker} \lambda_F$ is the quotient of a trivial bundle, namely $f^*(T_{Z/A}) \vert F$ and since $F$ is rationally connected,
this clearly contradicts $\det Q' = 0.$ So $\lambda_F = 0$ and hence $\lambda = 0.$ 
Therefore we obtain an epimorphism
$$ f^*(T_{Z/A} = p_B^*(T_B) \to Q'.$$
But $T_B$ is generically nef, and so does $p_B^*(T_B)$. This contradicts the negativity of $Q'$ on certain curves through $x.$
\end{proof}

\section{Threefolds} 
\setcounter{lemma}{0}

We end by considering Conjecture 1.3 in dimension 3. 

\begin{theorem} Let $X$ be a smooth projective threefold with $-K_X$ nef. Then $T_X$ is generically nef for some $(H_1,H_2) $ and
even generically nef for all $(H_1,H_2)$ 
unless (possibly) we are in the
the following case: \\
$X$ is rationally connected, $-K_X$ is ample on all movable curves; moreover there is a fibration $f: X \to B \simeq \bP_1,$ whose
general fiber has non-semi-ample anti-canonical bundle. 
\end{theorem}

\begin{proof}  By (6.4) we may assume that $-K_X$ is not semi-ample. Thus by [BP04] we are in one of the following cases.
\begin{enumerate} 
\item There is a finite \'etale cover $\tilde X \to X$ such that $\tilde X \simeq B \times S$ with $B$ an elliptic curve and $-K_S$ is not semi-ample.  
\item $X$ is rationally connected, $-K_X$ is ample on all movable curves; moreover there is a fibration $f: X \to B \simeq \bP_1,$ whose
general fiber has non-semi-ample anti-canonical bundle. 
\item $-K_X$ is ample on all movable curves and $X$ is a $\bP_1-$bundle over an abelian surface.
\end{enumerate} 
In case (3) $T_X$ is even nef amd there is nothing to prove. In case (1) clearly $q(\tilde X) \leq 2.$ 
If $q(\tilde X) = 2,$ so that $q(S) = 1,$ then the Albanese of $\tilde X$, being smooth by [PS98], realises $\tilde X$ is
a $\bP_1-$bundle over an abelian surface. Hence $T_{\tilde X}$ is nef and so does $T_X$. So we may assume $q(S) = 0$ and thus $S$ is $\bP_2$ blown up
in 9 points in general position. By (5.11) $T_S$ is generically ample, hence $T_{\tilde X}$ is easily seen to be generically nef. 
\vskip .2cm \noindent We finally treat the most difficult case (2) and need to prove the existence of some $h = (H_1,H_2)$ for which $T_X$ is generically nef. We
fix an ample divisor $H$ on $X$ and set $A = \sO_B(1).$ Introduce
$$ h_0 = p^*(A) \cdot H $$
and notice that 
$$K_X \cdot h_0 < 0.$$
We claim that there is an open neighborhood $U$ of $h_0$ in $\overline{ME}(X)$ such that 
for all line bundles $L \subset T_X$ and all $h$ with $[h] \in U,$
$$ (K_X + L) \cdot h < 0. \eqno (3.a)$$ 
For the proof of (3.b) we first observe that
$$ (K_X + L) \cdot h_0 < 0. \eqno (3.b)$$
In fact, by the definition of $h_0$ this comes down to 
$$ (K_F + L_F) \cdot H_F < 0. \eqno (*)$$ 
If the composed map $K_F \otimes L_F \to K_F \otimes T_X \vert F \to K_F \otimes N_F = K_F$ 
is non-zero, then (*) and (3.c) are already clear. If this map vanishes, then we have an inclusion 
$$ K_F \otimes L_F \subset K_F \otimes T_F = \Omega^1_F.$$
Since $-K_F$ is not semi-ample, $F$ is either the plane blown up in 9 points in sufficiently general position
or the projectivisation of a suitable semi-stable rank 2-bundle on an elliptic curve. In the first case $T_X$ is generically
ample, hence (*) follows. In the second case $T_F$ is nef and (*) is also easily checked. This establishes (3.c). \\
From (3.b) the claim (3.a) is an immediate consequence since the set
$$ \{ c_1(L) \vert  L \subset T_X; \ L \cdot h_0 > 0 \} $$
is finite. 
\vskip .2cm \noindent 
We claim that $T_X$ is $h$-generically nef for $[h] \in U(h_0).$ 
Supposing the contrary for a fixed $h,$ we find by (5.4) a subsheaf $\sE \subset T_X$ such that 
$$ (K_X + \det \sE) \cdot h > 0 \eqno (**)$$
and $\E$ is $h_{\epsilon}-$maximally nef. \\
But then (3.a) conflicts with (**), so we obtain a contradiction and conclude. 

\end{proof} 

\begin{theorem} Let $X$ be a smooth projective threefold with $-K_X$ nef. Then $T_X$ is suffciently nef.
\end{theorem}

\begin{proof} If $-K_X$ is hermitian semi-positive, we simply apply (7.11). Otherwise we are in of the three ``exceptional'' cases
listed in the proof of (8.1). In cases (1) and (3) it is clear that $T_X$ is sufficiently nef; in case (2) we can argue as in (8.1) using (7.8) 
and everything remains the same except that the positivity (**) is replaced
by $$ K_X + ßdet \sE \equiv 0.$$ 

\end{proof}


\begin{thebibliography}{Mum69}

\bibitem[BCHM06]{BCHM06} Birkar,C; Cascini,P.; Hacon,C.; McKernan,J.: 
Existence of minimal models for varieties of log general type.
arXiv:math/0610203 

\bibitem[BM01]{BM01} Bogomolov,F.; McQuillan,M.: Rational curves on foliated varieties. Preprint, IHES, 2001. 

\bibitem[BP04]{BP04} Bauer,Th. Peternell,Th.: 
Nef reduction and anticanonical bundles. 
Asian J. Math. 8, No. 2 (2004), 315--352

\bibitem[Bo04]{Bo04} Boucksom, S.: Divisorial Zariski decompositions on compact complex manifolds.  Ann. Sci. École Norm. Sup. (4)  37  (2004),  no. 1, 45--76. 

\bibitem[BP07]{BP07} Berndtsson,B.; Paun,M.: 
Bergman kernels and the pseudoeffectivity of relative canonical bundles.  arXiv:math/0703344

\bibitem[CP98]{CP98} Campana,F.; Peternell,Th.: Rational curves and ampleness properties of the tangent bundle of algebraic varieties.  
Manuscripta Math.  97  (1998),  no. 1, 59--74. 
 

\bibitem[CP07]{CP07} Campana,F.; Peternell, Th.: Geometric stability of the cotangent bundle and the universal cover of a projective manifold
(with an appendix by Matei Toma).   arXiv:math/0405093 

\bibitem[DPS94]{DPS94}  Demailly,J.P.; Peternell,Th.; Schneider,M.:
Compact complex manifolds with numerically effective tangent bundles. 
J. Algebr. Geom. 3, No.2, 295-345 (1994).


\bibitem[DPS96]{DPS96} Demailly,J.P.; Peternell,Th.; Schneider,M.: Compact K\"ahler manifolds with Hermitian semipositive anticanonical bundle.  
Compositio Math.  101  (1996),  no. 2, 217--224.

\bibitem[DPS01]{DPS01} Demailly, J.P.; Peternell, Th.; Schneider, M.: 
Pseudo-effective line bundles on compact Kähler manifolds. 
Int. J. Math. 12, No. 6, 689-741 (2001). 

\bibitem[Fl88]{Fl88} Flenner,H.: Extendability of differential form on non-isolated singularities. Inv.math. 94 (1988), 317--326

\bibitem[Fu86]{Fu86} Fujita,T.: Zariski decomposition and canonical rings of elliptic threefolds. 
J. Math. Soc. Japan, 38, No.1, (1986), 19--37 

\bibitem[Hoe06]{Hoe06} H\"oring, Andreas: Two applications of positivity to the classification theory of complex projective manifolds. 
Ph. D. thesis, Bayreuth, Grenoble 2006 

\bibitem[Ka92]{Ka92} Kawamata, Y.: Abundance theorem for minimal threefolds.  Invent. Math.  108  (1992),  no. 2, 229--246. 

\bibitem[KST07]{KST07} Kebekus,S.; Sola Cond\'e,L; Toma,M.: Rationally connected foliations after Bogomolov and McQuillan.  
J. Algebraic Geom.  16  (2007),  no. 1, 65--81. 

\bibitem[MM86]{MM86} Miyaoka,Y.; Mori,S.: A numerical criterion for uniruledness.  Ann. of Math. (2)  124  (1986),  no. 1, 65--69.

\bibitem[Mi87]{Mi87} Miyaoka,Y.: The Chern classes and Kodaira dimension of a minimal variety.  Algebraic geometry, Sendai, 1985,  449--476, 
Adv. Stud. Pure Math., 10, North-Holland, Amsterdam, 1987.

\bibitem[MR82]{MR82} Mehta, V. B.; Ramanathan, A. Semistable sheaves on projective varieties and their restriction to curves.  
Math. Ann.  258  (1981/82), no. 3, 213--224. 

\bibitem[Pe94]{Pe94} Peternell, Th.
Minimal varieties with trivial canonical classes. I. 
Math. Z. 217, No.3, 377-405 (1994). 

\bibitem[Pe06]{Pe06} Peternell, Th.:
Kodaira dimension of subvarieties. II. 
Int. J. Math. 17, No. 5, 619-631 (2006).

\bibitem[PS98]{PS98} Peternell,Th.; Serrano, F.:
Threefolds with nef anticanonical bundles. 
Collect. Math. 49, No.2-3, 465-517 (1998). 

\bibitem[PS00]{PS00} Peternell, Th.; Sommese, S.J.:  Ample vector bundles and branched coverings. Comm. Algebra  28  (2000),  no. 12, 5573--5599. 

\bibitem[RC00]{RC00} Ran, Z.; Clemens, H.: A new method in Fano geometry.  Internat. Math. Res. Notices  (2000),  no. 10, 527--549.

\bibitem[Si74]{Si74} Siu,Y.T.: Analyticity of sets associated to Lelong numbers and the extension of closed positive currents.  Invent. Math.  27  (1974), 53--156. 

\bibitem[Zh05]{Zh05} Zhang, Qi: On projective manifolds with nef anticanonical divisors. Math. Ann. 332 (2005), 697--703

\bibitem[Zh06]{Zh06} Zhang, Qi: Rational connectedness of log Q-Fano varieties. J. reine u. angew. Math.  590  (2006), 131--142.


\end{thebibliography}
\end{document}